\documentclass[runningheads]{llncs}

\usepackage[utf8]{inputenc} 
\usepackage[T1]{fontenc}    
\usepackage{hyperref}       
\usepackage{url}            
\usepackage{booktabs}       
\usepackage{amsfonts}       
\usepackage{nicefrac}       
\usepackage{microtype}      
\usepackage{xcolor}         
\usepackage{amsmath,amsfonts,amssymb,array}
\usepackage{algorithm}
\usepackage[font=small]{caption}
\usepackage{subcaption}
\usepackage{algpseudocode}
\usepackage{bm}
\usepackage{color}
\usepackage{graphicx}
\usepackage{enumitem}

\newtheorem{assumption}{Assumption}

\newcommand{\R}{\ensuremath{\mathbb{R}}}
\newcommand{\N}{\ensuremath{\mathbb{N}}}

\newcommand{\hook}{\ensuremath{\hookrightarrow}}

\begin{document}

\title{Extragradient Sliding for Composite Non-Monotone Variational Inequalities}
\titlerunning{Extragradient Sliding for Composite Non-Monotone VIs}
%

\author{Roman Emelyanov \inst{1} \and
Andrey Tikhomirov\inst{1} \and
Aleksandr Beznosikov\inst{1,2,3}
\and
Alexander Gasnikov\inst{3,1,2}
}

\authorrunning{R. Emelyanov, A. Tikhomirov, A. Beznosikov, A. Gasnikov}
%
\institute{Moscow Institute of Physics and Technology, Moscow, Russian Federation \and
Institute for Information Transmission Problems RAS, Moscow, Russian Federation
 \and
Innopolis University, Innopolis, Russian Federation}
\maketitle              
\begin{abstract}
Variational inequalities offer a versatile and straightforward approach to analyzing a broad range of equilibrium problems in both theoretical and practical fields. 
In this paper, we consider a composite generally non-monotone variational inequality represented as a sum of $L_q$-Lipschitz monotone and $L_p$-Lipschitz generally non-monotone operators. We applied a special sliding version of the classical Extragradient method to this problem and obtain better convergence results. In particular, to achieve $\varepsilon$-accuracy of the solution, the oracle complexity of the non-monotone operator $Q$ for our algorithm is $\mathcal{O}\left(L_p^2/\varepsilon^2\right)$ in contrast to the basic Extragradient algorithm with $\mathcal{O}\left((L_p+L_q)^2/\varepsilon^2\right)$.
The results of numerical experiments confirm the theoretical findings and show the superiority of the proposed method.

\end{abstract}

\section{Introduction}\label{sec:intro}


Variational inequalities (VIs) are a popular class of optimization problems, which despite its relative youth has an extensive history of research, both in terms of different formulations, and of effective methods of their solution. The variational inequalities paradigm has gained particular popularity due to its generality and its ability to describe and represent various optimization problems in a unified way \cite{HarkerVIsurvey1990,VIbook2003}. Nowadays VI problems can be found in a large number of fields from economics, game theory, physics and modelling of transport flows to machine learning and rapidly developing deep learning \cite{jofre2007variational,facchinei2003finite,shafieezadeh2015distributionally}. For instance, the development of optimization methods aimed at solving variational inequalities has recently attracted the attention of many researchers in the context of optimizing loss functions of generative adversarial networks \cite{goodfellow2014generative} and reinforcement learning \cite{pinto2017robust}. 

The most straightforward and widely known method for dealing with problems posed as variational inequalities is Gradient Descent-Ascent \cite{browder1966existence,rockafellar1969convex,sibony1970methodes} developed by analogy with methods of optimization of a single objective \cite{polyak1987introduction}. Its serious drawbacks \cite{beznosikov2023smooth}, which include weaker convergence compared to ordinary gradient descent and the problem of rotation around the optimum \cite{polyak1987introduction,beznosikov2023smooth} led researchers to create more advanced methods.
One of the most well-known of these algorithms is Extragradient, proposed by G. Korpelevich \cite{korpelevich1976extragradient}. Over time, the research around this method evolved and there are now various modifications (e.g. Optimistic Gradient with one oracle call per iteration \cite{popov1980modification,mokhtari2020convergence}) and generalizations to an arbitrary Bregman setup (e.g. Mirror-Prox \cite{nemirovski2004prox}). Moreover, for the monotone variational inequalities, the optimality of the Extragradient method was also shown. Meanwhile, there is a large body of theoretical studies related to various kinds of relaxations of monotonicity and transition to non-monotone operators \cite{dang2015convergence,kannan2019optimal,diakonikolas2021efficient,beznosikov2022decentralized,beznosikov2022distributed1}. The aim of this paper is also to delve into the non-monotone setting, but to look deeper into it and consider composite operators, i.e., operators that are represented as the sum of two operators. It seems that using the additional structure of the target problem can give improvements in terms of convergence theory. For monotone and strongly monotone operators such approaches already exist and indeed give theoretical and practical improvements \cite{lan2021mirror,beznosikov2021distributed,kovalev2022optimal}.



\paragraph{Our contribution.} We consider the variational inequality with a composite operator $R := P + Q$, both components of which are Lipschitz-continuous and only one of the component $Q$ is monotone (the operator $P$ can be generally non-monotone). We additionally assume that the whole target operator $R$ satisfies Minty assumption \cite{minty62} -- by far the most common and well-known relaxation of non-monotonicity. For this kind of problem, the classical Extragradient method requires \cite{dang2015convergence} 
$$
\mathcal{O}\left( \frac{(L_p + L_q)^2 \| x^0 - x^*\|^2}{\varepsilon^2} \right) \text{computations of the operator } P,
$$
where $L_p$ and $L_q$ are Lipschitz constants of $P$ and $Q$, $x^0$ is a starting point, $x^*$ is a solution of the VI problem, $\varepsilon$ is required accuracy of the obtained numerical solution. On the other hand, motivated by the fact that the non-monotone operator $P$ is likely to be more computationally expensive than $Q$, we consider a modification of Extragradient using the sliding technique \cite{kovalev2022optimal}, which allows to take into account the composite structure of the problem. This makes it possible to compute one of the operators less frequently. In particular, we obtain improvements on the estimate of operator $P$ calls. For our method it is
$$
\mathcal{O}\left( \frac{L_p^2 \| x^0 - x^*\|^2}{\varepsilon^2} \right),
$$
which can be much better for some relations on $L_p$ and $L_q$ compared to the Extragradient method. To support the theoretical results, a series of experiments are set up and the results confirm the improvements.

\paragraph{Notation.} We use $\langle x,y \rangle := \sum_{i=1}^d x_i y_i$ to introduce inner product of $x,y\in\R^d$, where $x_i$ corresponds to the $i$-th component of $x$ in the standard basis in $\R^d$. It induces $\ell_2$-norm in $\R^d$ according to $\|x\| := \sqrt{\langle x,x \rangle}$.

\section{Problem setting}
We consider the composite variational inequality in the following form: 
\begin{equation}
    \label{eqn::problem}
    \text{Find } ~~x^* \in \R^d : R(x^*) = 0 ~~\text{ with }~~ R(x) := Q(x) + P(x),
\end{equation}
where $Q(x), P(x): \mathbb{R}^d \rightarrow \mathbb{R}^d$ are operators.

Many known problems can be reformulated using the language of composite variational inequalities. Let us consider two common use cases for VIs:
\begin{enumerate} 
    
    \item {\bf Minimization problem.}
    One can notice, that solving
    $\min_{x\in \mathbb{R}^d} r(x)$,
    where $r(x) := q(x)+p(x)$ is a convex function,
    is equivalent to solving the VI problem \eqref{eqn::problem} with $Q(x):=\nabla q(x)$, $P(x) := \nabla p(x)$, $R(x) := \nabla r(x) = \nabla q(x) + \nabla p(x)$.

    \item {\bf Saddle point problem.}
    Let us consider the saddle-point problem:
    \begin{equation}
        \label{eqn::SSP}
        \min_{y \in \R^{d_y}} \max_{z \in \R^{d_z}} [r(y,z) := q(y,z) + p(y,z)].
    \end{equation}

    If we take $x = [y,z]$, $Q(x) := Q(y, z) = [\nabla_y q(y,z), -\nabla_z q(y,z)]$, $P(x) := P(y, z) = [\nabla_y p(y,z), -\nabla_z p(y,z)]$, then it can be proved for the convex-concave function $r(y,z)$ that $x^* = (y^*, z^*)$ is a solution for \eqref{eqn::problem} if and only if the following inequality holds
    $$
    r(y^*, z) \leq r(y^*, z^*) \leq r(y, z^*) \; \forall y \in \R^{d_y}, z \in \R^{d_z}.
    $$
    Equivalently, it means that $x^* = [y^*, z^*]$ is a solution for \eqref{eqn::SSP}.
    
\end{enumerate}

Despite the fact that a minimization problem is a special case of variational inequalities, they are usually studied separately. This is due to the fact that a more optimistic convergence theory can be constructed for minimization problems \cite{polyak1987introduction,nesterov2018lectures} compared to the general results for variational inequalities. 
Meanwhile, the study of saddle point problems is often conducted through the prism of variational inequalities \cite{nemirovski2004prox,juditsky2011solving}, particularly in recent years, theoretical studies of variational inequalities have been associated with solving the practical minimax learning problem of GANs training \cite{gidel2018variational,beznosikov2022decentralized,beznosikov2022scaled}.

We study problem \eqref{eqn::problem} under the following commonly used assumptions:

\begin{assumption}
\label{asm::R}
$R(x)$ satisfies Minty assumption:
$$ \exists\: x^* \in \R^d: \;  \forall x \in \R^d \hook \langle R(x), x - x^* \rangle \geq 0.$$
\end{assumption}

This assumption, also called the variational stability condition is considered as an option to structurally constrain a non-monotone problem. Minty assumption is widely used in the literature \cite{dang2015convergence,iusem2017extragradient,mertikopoulos2018optimistic,hsieh2020explore,diakonikolas2021efficient,beznosikov2022distributed}.

Next, we also introduce two standard assumptions for the analysis of variational inequalities.

\begin{assumption}
\label{asm::Q}
$Q(x)$ is $L_q$-Lipschitz and monotone:
$$\forall\: x_1, x_2 \in \R^d \hook \|Q(x_1) - Q(x_2)\| \leq L_q \|x_1 - x_2\|, $$
$$ \forall\: x_1, x_2 \in \R^d \hook \langle Q(x_1) - Q(x_2), x_1 - x_2 \rangle \geq 0. $$
\end{assumption}

\begin{assumption}
\label{asm::P}
$P(x)$ is $L_p$-Lipschitz:
$$ \forall\: x_1, x_2 \in \R^d \hook \|P(x_1) - P(x_2)\| \leq L_p \|x_1-x_2\|. $$
\end{assumption}

Once again we emphasize the key detail that only the operator $Q$, but not $P$, is monotone, and hence the full operator $R$ can be non-monotone in general. 

\section{Algorithm}

The algorithm studied in this paper is a version of the Extragradient algorithm, but with additional sliding technique \cite{kovalev2022optimal}:

\begin{algorithm}[!h]
    \caption{Extragradient Sliding} \label{alg1}
    \begin{algorithmic}[1]
        \State $\textbf{Input: }$ starting point $x^0 \in \R^d$
        \State $\textbf{Parameters: }$ stepsizes $ \eta, \theta > 0$, number iterations $K \in \N $
        \For{$k = 0, 1, 2, \dots, K-1$}
            \State Find $u^k \approx \tilde u^k$ where $\tilde u^k$ is solution for 
            \Statex
            \hspace{1cm} Find $\tilde u^k \in \R^d: B_{\theta}^k(\tilde u^k) = 0$ with $B_{\theta}^k(x) := P(x^k) + Q(x) + \frac{1}{\theta}(x - x^k)$ \label{alg::subproblem}
            \State $x^{k+1} = x^k - \eta R(u^k)$ \label{alg::step}
        \EndFor 
    \end{algorithmic}
\end{algorithm}

Let us give a high-level intuition of how the above algorithm works. The main idea of this algorithm is to move away from equal number of calls of $P$ and $Q$, as it happens in the classical Extragradient method. The more computationally expensive $P$ is called twice per iteration of the algorithm, when selecting the optimal $\tilde u^k$ in line \ref{alg::subproblem} and when computing $R$ in line \ref{alg::step}. In turn, the computationally simpler $Q$ is called some number of times in the inner problem (line \ref{alg::subproblem}) and also when computing $R$. This is the idea behind the sliding technique: fix one of the operators and vary the other due to its cheapness. Thus in line \ref{alg::subproblem} the full operator $R(x)$ is approximated by $Q(x)$ and a slightly outdated version $P(x^k)$. 



\section{Convergence analysis}

\begin{theorem}
\label{th::main}
Consider Algorithm \ref{alg1} for the problem \eqref{eqn::problem} under Assumptions \ref{asm::R}--\ref{asm::P}, with the following tuning:
\begin{equation*}
  \theta = \dfrac{1}{2L_p},\; \eta = \dfrac{\theta}{2}.  
\end{equation*}
Assume that $u^k$ (line \ref{alg::subproblem}) satisfies:
\begin{equation}
\label{asm:B}
     \|B_{\theta}^k(u^k) \|^2 \leq \dfrac{L_p^2}{3} \|x^k - \tilde u^k \|^2.
\end{equation}
Then, we have the following convergence estimate:
\begin{equation*} 
     \underset{0 \leq j \leq K-1}{\min} \|R(u^j) \|^2 \leq \dfrac{16L_p^2 \|x_0 - x^*\|^2}{K}.
\end{equation*}
\end{theorem}
This results means sublinear convergence. To prove the theorem we first deal with the auxiliary lemma.

\begin{lemma}
\label{lem::1}
Consider Algorithm \ref{alg1}. Let $\theta$ be defined as $\theta = \tfrac{1}{2L_p}$. Then, under Assumptions \ref{asm::R}, \ref{asm::Q}, \ref{asm::P}, the following inequality holds:
\begin{equation}
    2\langle x^* - x^k, R(u^k) \rangle \leq - \theta \|R(u^k)\|^2 + 3\theta \left( \|B_{\theta}^k(u^k)\|^2 - \frac{L_p^2}{3}\|x^k- \tilde u^k \|^2 \right).
\end{equation}
\end{lemma}

\noindent
{\bf Proof of Lemma \ref{lem::1}.}
Using Assumption \ref{asm::R}, we get
\begin{eqnarray*}
    2\langle x^k - x^*, R(u^k)\rangle &=& 2\langle x^* - u^k, R(u^k) \rangle + 2\langle u^k - x^k, R(u^k) \rangle
    \\ &\leq& 
    2\langle u^k - x^k, R(u^k) \rangle = 2\theta \left \langle \dfrac{1}{\theta} (u^k - x^k), R(u^k) \right \rangle.
\end{eqnarray*}
The definition of $B_{\theta}^k(x)$ (line \ref{alg::subproblem} of Algorithm \ref{alg1}) gives
\begin{eqnarray*}
    2\theta \left \langle \dfrac{1}{\theta} (u^k - x^k), R(u^k) \right \rangle
    &=&
    \theta \left \| \dfrac{1}{\theta} (u^k - x^k) + R(u^k) \right\|^2 - \dfrac{1}{\theta} \|u^k - x^k \|^2 - \theta \| R(u^k) \|^2 
    \\&=&
    -\dfrac{1}{\theta} \|u^k - x^k \|^2 - \theta \|R(u^k) \|^2 
    \\
    &&+ \theta \|B_{\theta}^k(u^k) - P(x^k) + P(u^k)\|^2.
\end{eqnarray*}
Using the Cauchy-Schwarz inequality and $L_p$-Lipschitzness of $P(x)$ (Assumption \ref{asm::P}), we get
\begin{eqnarray*}
    2\theta \left \langle \dfrac{1}{\theta} (u^k - x^k), R(u^k) \right \rangle
    &\leq&
    -\dfrac{1}{\theta} \|u^k-x^k\|-\theta \|R(u^k)\|^2 + 2\theta \|B_{\theta}^k(u^k)\|^2 
    \\
    &&+ 2 \|P(u^k) - P(x^k) \|^2
    \\&\leq&
    -\dfrac{1}{\theta} \|u^k-x^k\|-\theta \|R(u^k)\|^2 + 2\theta \|B_{\theta}^k(u^k)\|^2 
    \\
    &&+ 2\theta L_p^2 \|u^k - x^k\|^2
    \\&=&
    - \dfrac{1}{\theta} \left( 1-2\theta^2 L_p^2 \right)\|u^k-x^k\|^2 - \theta \| R(u^k)\|^2 + 2\theta \|B_{\theta}^k(u^k)\|^2.
\end{eqnarray*}
With $\theta = \tfrac{1}{2L_p}$ and the Cauchy-Schwarz inequality in the form $- \|a\|^2 \leq \|b\|^2 - \dfrac{1}{2} \|a + b \|^2$, one can obtain
\begin{eqnarray*}
    2\langle x^* - x^k, R(u^k) \rangle 
    &\leq&
    - \theta \|R(u^k)\|^2 + 2\theta \|B_{\theta}^k(u^k)\|^2 -\dfrac{1}{2\theta}\|u^k-x^k\|^2
    \\&\leq &
    -\theta\|R(u^k)\|^2 + 2\theta \|B_{\theta}^k(u^k)\|^2 
    \\
    &&+ \dfrac{1}{2\theta}\|u^k - \tilde u^k\|^2 - \dfrac{1}{4\theta}\|x^k - \tilde u^k \|^2.
\end{eqnarray*}

Additionally, we can observe that $B_\theta^k (x)$ is $\tfrac{1}{\theta}$-strongly monotone. It follows from the definition of the operator: the operator $B_\theta^k (x)$ is a sum of the monotone operator $Q$ (Assumption \ref{asm::Q}) and the strong monotone linear operator  $\frac{1}{\theta}(x - x^k)$. Together with the Cauchy-Schwarz inequality, it gives that
 \begin{eqnarray*}
     \|u^k - \tilde u^k\|^2 \leq
     \theta \langle B_{\theta}^k(u^k) - B_{\theta}^k(\tilde u^k), u^k - \tilde u^k \rangle \leq
     \theta \| B_{\theta}^k(u^k) - B_{\theta}^k(\tilde u^k) \| \cdot \|u^k - \tilde u^k\|.
 \end{eqnarray*} 
With $B_{\theta}^k(\tilde u^k) = 0$ ($\tilde u^k$ is the solution of the subproblem from line \ref{alg::subproblem}), we get
 \begin{eqnarray*}
     \|u^k -\tilde u^k\|^2 \leq \theta^2 \|B_{\theta}^k(u^k)\|^2.
 \end{eqnarray*}
Applying this to the upper inequality, we finalize the proof:
\begin{eqnarray*}
    2\langle x^* - x^k, R(u^k) \rangle 
    &\leq&
    - \theta \|R(u^k)\|^2 + \dfrac{5}{2}\theta \|B_\theta^k(u^k)\|^2 
    - \dfrac{1}{4\theta}\|x^k - \tilde u^k \|^2 
    \\&\leq& 
    - \theta \|R(u^k)\|^2 + 3\theta \|B_{\theta}^k (u^k) \|^2 - \dfrac{3 \theta}{12 \theta^2} \|x^k - \tilde u^k \|^2 
    \\&=&
    - \theta \|R(u^k)\|^2 + 3\theta \left(\|B_{\theta}^k (u^k) \|^2 - 
    \dfrac{L_p^2}{3} \|x^k - \tilde u^k \|^2  \right).
\end{eqnarray*}
Here we substitute $\theta = \tfrac{1}{2L_p}$. $\square$

Now we are ready to prove the main theorem.

\noindent
{\bf Proof of Theorem \ref{th::main}.}
Line \ref{alg::step} of Algorithm \ref{alg1} gives
\begin{eqnarray*}
    \|x^{k+1}-x^*\|^2
    &=&
    \|x^{k+1}-x^k\|^2+2\langle x^{k+1}-x^k, x^k - x^* \rangle + \|x^k - x^*\|^2
    \\&=&
    \|x^{k+1}-x^k\|^2 + \|x^{k}-x^*\|^2  - 2 \eta \langle R(u^k), x^k - x^* \rangle.
\end{eqnarray*}
Using the results of Lemma \ref{lem::1} and the condition \eqref{asm:B} on $\|B_{\theta}^k(u^k)\|^2$, we get 
\begin{eqnarray*}
    \|x^{k+1}-x^*\|^2
    &\leq&
    \|x^{k+1}-x^k\|^2 + \|x^k-x^*\|^2 - \eta \theta\|R(u^k)\|^2
    \\
    &&+
    3\eta \theta\left( \|B_{\theta}^k(u^k)\|^2-\dfrac{L_p^2}{3}\|x^k-\tilde u^k\|^2 \right)
    \\
    &\leq&
    \|x^{k+1}-x^k\|^2 + \|x^k-x^*\|^2 - \eta \theta \|R(u^k)\|^2.
\end{eqnarray*}
Again from line \ref{alg::step} it follows that
\begin{eqnarray*}
    \|x^{k+1}-x^*\|^2
    &\leq& 
    \eta^2 \|R(u^k)\|^2 + \|x^k - x^*\|^2 - \eta \theta \|R(u^k)\|^2.
\end{eqnarray*}
Let us substitute the choice of parameters: $\theta=\tfrac{1}{2L_p}$, $\eta = \tfrac{\theta}{2}$:
\begin{eqnarray*}
    \|x^{k+1}-x^*\|^2 &\leq&
    \|x^k-x^*\|^2 - \eta (\theta - \eta) \|R(u^k)\|^2
    \\&\leq&
    \|x^k-x^*\|^2 - \dfrac{\theta^2}{4}\|R(u^k)\|^2.
\end{eqnarray*}
Summing from 1 to $K - 1$, one can obtain
\begin{eqnarray*}
    \sum\limits_{j=0}^{K-1} \dfrac{\theta^2}{4}\|R(u^j)\|^2 
    &\leq&
    \sum \limits_{j=0}^{K-1} \left( \|x^j-x^*\|^2-\|x^{j+1}-x^*\|^2 \right)
    \\&=& 
    \|x^0-x^*\|^2-\|x^K-x^*\|^2 
    \\
    &\leq& \|x^0 - x^*\|^2.
\end{eqnarray*}
Thus, we have
$$\sum\limits_{j=0}^{K-1}\|R(u^j)\|^2\leq 16L_p^2\|x^0-x^*\|^2,$$
and finally,
$$\underset{0 \leq j \leq K-1}{\min} \|R(u^j)\|^2 \leq \dfrac{16L_p^2\|x^0-x^*\|^2}{K}.$$
This ends the proof. $\square$

\begin{corollary}
Under the assumptions of Theorem \ref{th::main}, to achieve a $\varepsilon$-solution in terms $\varepsilon \sim \| R(u)\|$, Algorithm \ref{alg1} needs 
$$
\mathcal{O}\left( \frac{L_p^2 \| x^0 - x^*\|^2}{\varepsilon^2}\right) \text{ computations of the operator } P.
$$
\end{corollary}
As noted in Section \ref{sec:intro}, this result is better than for the Extragradient method, which provides an estimate: $
\mathcal{O}\left( \tfrac{(L_p+L_q)^2 \| x^0 - x^*\|^2}{\varepsilon^2}\right).
$






\begin{remark} 
Meanwhile, line \ref{alg::subproblem} of the algorithm requires an additional algorithm to efficiently solve the resulting subproblem. The Extra Anchored Gradient algorithm proposed in \cite{yoon2021accelerated} can be used for these purposes. In fact, any method of solving variational inequalities with a Lipschitz monotone operator can be used here, but the method from \cite{yoon2021accelerated} has convergence guarantees necessary for our theoretical analysis (see \eqref{asm:B}). Moreover, it is also shown in \cite{yoon2021accelerated} that these guarantees are optimal and unimprovable.
\end{remark}

\section{Numerical Experiments}

In this part, we conduct three experiments: with a generated bilinear problem, with a logarithmic logistic regression, and with non-convex least squares (NLLSQ) on the \texttt{mushrooms} dataset from LibSVM \cite{CC01a,misc_mushroom_73}. It turns out we evaluate the performance of the algorithm on both synthetic and real data.

\subsection{Bilinear problem}
A bilinear problem is a classical and keystone example of the saddle:
\begin{equation}
\label{eq:bilinear}
\min_{x \in [-1; 1]^d} \max_{y \in [-1; 1]^d} \left[f(x,y) := (x-b_x)^T A (y - b_y) + \frac{1}{2}\| x - b_x\|^2 - \frac{1}{2}\| y - b_y\|^2 \right].
\end{equation}
In this setting, $P$ refers to the main part $(x-b_x)^T A (y - b_y)$ of the gradient while $Q$ represents the gradient from the regularization terms $\tfrac{1}{2}\| x - b_x\|^2 - \tfrac{1}{2}\| y - b_y\|^2$ (see the second example from Section 2).

\begin{figure}[h]
\centering
\caption{Comparison of Extragradient and Extragradient Sliding for the generated bilinear saddle point problem \eqref{eq:bilinear}.}
\label{fig:comp1}
\begin{minipage}[][][b]{\textwidth}
\centering
\includegraphics[width=0.32\textwidth]{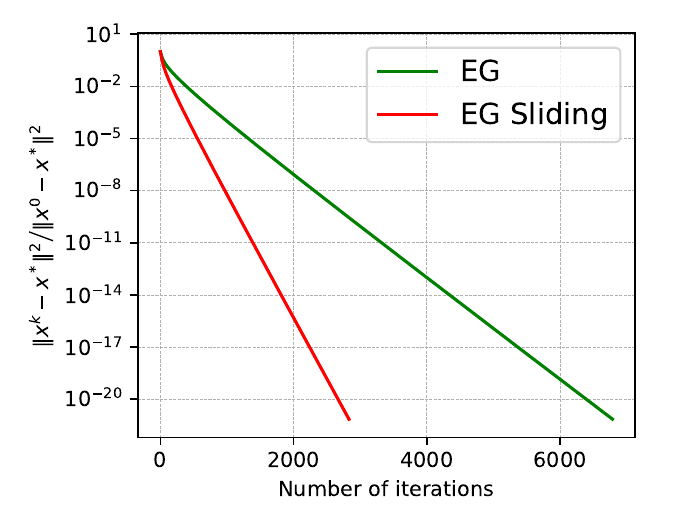}
\includegraphics[width=0.32\textwidth]{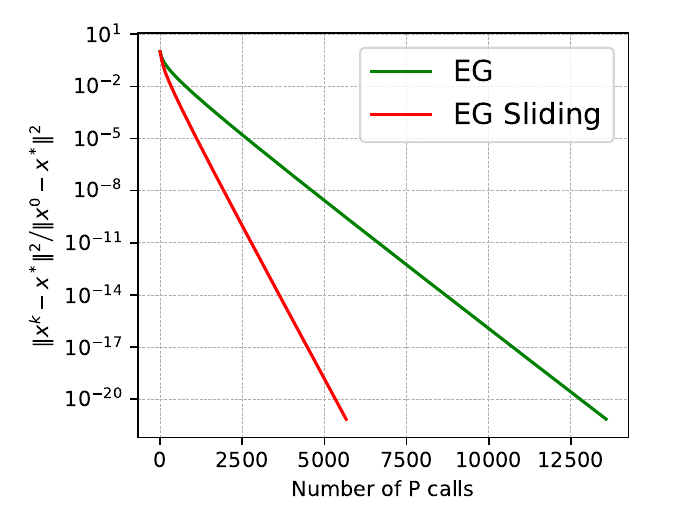}
\includegraphics[width=0.32\textwidth]{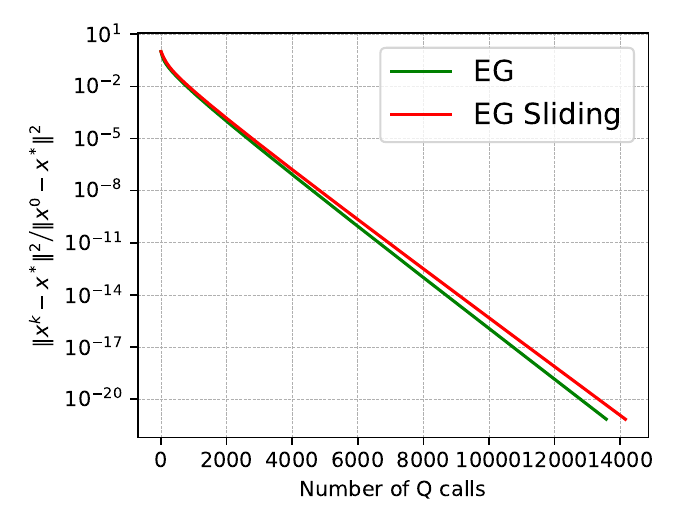}
\end{minipage}
\end{figure}

For the purposes of the experiment, a random bilinear saddle point problem is generated. The dimension $d$ of the problem is set equal to $1000$. The matrix $A$ are sampled as a positive definite matrix with uniformly distributed eigenvalues from $\mu$ to $L$, where $\mu$ and $L$ are chosen as $0.1$ and $100$ correspondingly. Biases $b_x$ and $b_y$ are both sampled from $\mathcal{U}(-1,1)$ as well as the starting point. 
In Figure \ref{fig:comp1}, one can see the plots comparing the convergence of the algorithm presented in this paper and Extragradient in terms of the number of iterations and oracle calls. 

\subsection{Logarithmic loss problem}

The second experiment uses the \textit{mushrooms} dataset and a more complex to compute logarithmic regression problem. To create a saddle point problem, the adversarial noise \cite{zhou2022modeling} technique is used. In this setting, the model weights are adjusted in parallel with the trained noise, which makes the model more robust. The resulting problem is formulated as follows:

\begin{equation}
\label{eq:logloss}
\begin{split}
    \min_{x \in \R^d} \max_{\|y_i\|\leq \delta} \bigg[f(x, y_1, \ldots y_i,\ldots y_N) :=& \dfrac{1}{N} \sum \limits_{i = 1}^{N} \ln \left(1 + \exp (-b_i x^T (A_i + y_i)) \right) 
\\
&+ \dfrac{\beta_x}{2} \|x\|^2 - \dfrac{\beta_y}{2} \|y\|^2 \bigg].
\end{split}
\end{equation}
Here $A, b$ are data, $x$ represents the model's weights, $y_i$ stands for adversarial noise, $\beta_x$ and $\beta_y$ determine the degree of regularization, $\delta$ defines the constraint imposed on adversarial noise. As in the previous case, the starting point is sampled from a uniform distribution $\mathcal{U}(-1,1)$. $\beta_x$ and $\beta_y$ are both set to $0.1$, $\delta$ is set to $0.1$. 

\begin{figure}[h]
\centering
\caption{Comparison of Extragradient and Extragradient Sliding for the log loss saddle point problem \eqref{eq:logloss}}
\label{fig:comp2}
\begin{minipage}[][][b]{\textwidth}
\centering
\includegraphics[width=0.32\textwidth]{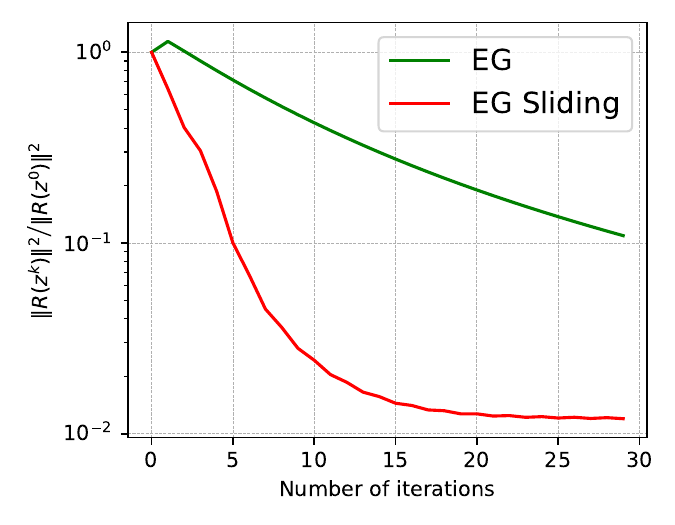}
\includegraphics[width=0.32\textwidth]{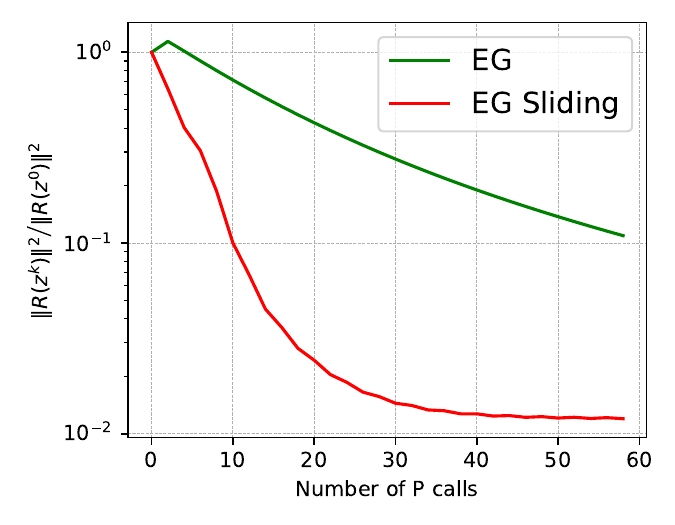}
\includegraphics[width=0.32\textwidth]{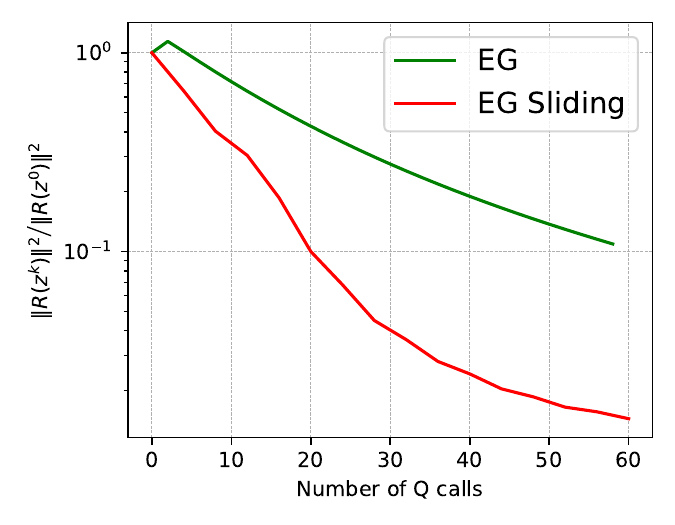}
\end{minipage}
\end{figure}

Plots showing comparison of the convergence for Algorithm \ref{alg1} and Extragradient are presented in Figure \ref{fig:comp2}. As can be seen, the algorithm presented in the paper outperforms baseline simultaneously in terms of $P$ and $Q$ oracle calls, despite the fact that solving the inner subproblem involves more $Q$ oracle calls than in the standard Extragradient.

\subsection{NLLSQ loss problem}

The last experiment is performed on a non-convex loss function in order to demonstrate the success of the method for the non-monotone operator $P$. As in the case of the logistic loss function, a saddle problem with adversarial noise adding robustness is posed:
\begin{equation}
\label{eq:NLLSQ}
\begin{split}
    \min_{x \in \mathbb{R}^d} \max_{\|y_i\|\leq \delta} \bigg[f(x,y_1, \ldots y_i,\ldots y_N) := & \dfrac{1}{N} \sum \limits_{i = 1}^{N} \left (b_i - \dfrac{1}{1+\exp (-x^T(A_i + y_i))}\right)^2 
    \\
    &+  \dfrac{\beta_x}{2} \|x\|^2 - \dfrac{\beta_y}{2} \|y\|^2 \bigg].
\end{split}
\end{equation}
The \textit{mushrooms} dataset is used again for the experiment. The notation used in the formula is similar to the previous case. The starting point is sampled from a uniform distribution $\mathcal{U}(-1,1)$. $\beta_x$ and $\beta_y$ are both set to $0.1$, $\delta$ is set to $0.1$.

\begin{figure}[h]
\centering
\caption{Comparison of Extragradient and Extragradient Sliding for the NLLSQ loss saddle point problem \eqref{eq:NLLSQ}}
\label{fig:comp3}
\begin{minipage}[][][b]{\textwidth}
\centering
\includegraphics[width=0.32\textwidth]{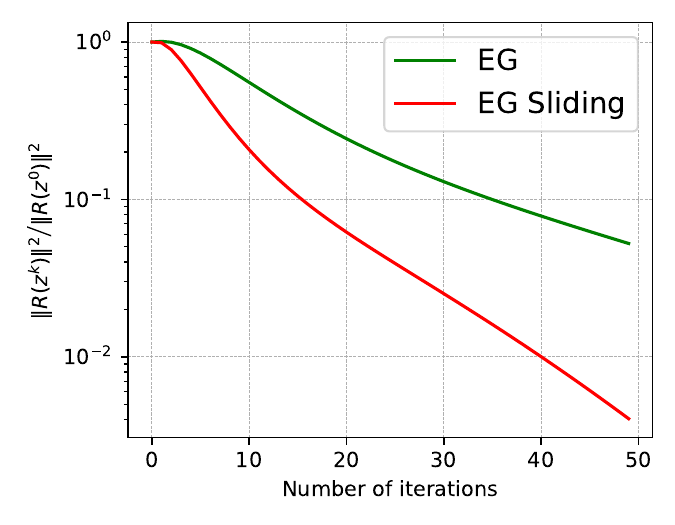}
\includegraphics[width=0.32\textwidth]{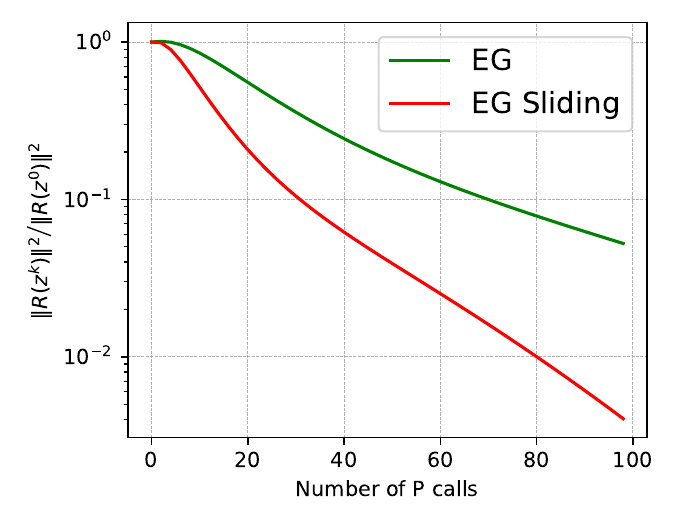}
\includegraphics[width=0.32\textwidth]{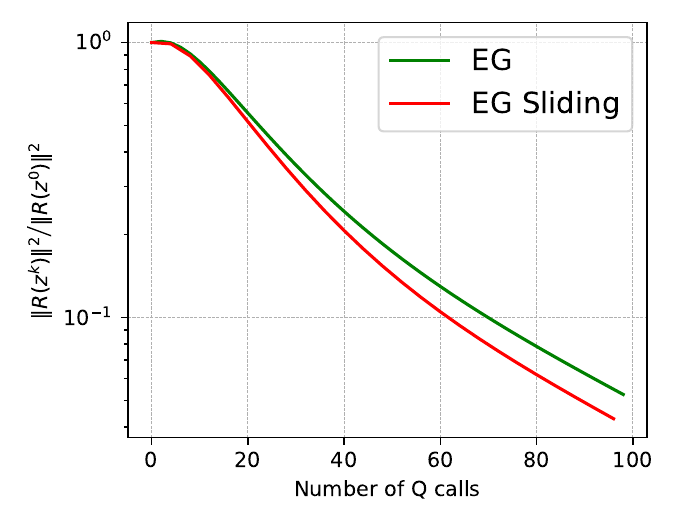}
\end{minipage}
\end{figure}

The comparison of algorithms in this setting, presented in Figure \ref{fig:comp3}, emphasizes the practical value of the considered algorithm. In this experiment it again wins the baseline on the calls of both oracles.



\bibliographystyle{unsrt}
\bibliography{refs}  
\end{document}